\documentclass[12pt,oneside,english]{amsart}
\usepackage[latin1]{inputenc}
\usepackage{amssymb}

\makeatletter
\newtheorem{theorem}{Theorem}

\usepackage{babel}

\makeatother
\begin{document}
\baselineskip=17pt

\title{Two digit theorems}

\author{Vladimir Shevelev}
\address{Departments of Mathematics \\Ben-Gurion University of the
 Negev\\Beer-Sheva 84105, Israel. e-mail:shevelev@bgu.ac.il}


\begin{abstract}
We prove that if $p$ is a prime with a primitive root 2 then
$S_{p,0}(2^p)=p$ and give a sufficient condition for an equality of
kind $S_{p,0}(2^p)=\pm p$.
\end{abstract}

\maketitle
\section{Introduction}
Denote for $x,m, l\in\mathbb{N},\;\;l\in[0,m-1]$,

$$
S_m(x)=S_{m,0}(x)=\sum_{0\leq n <x,\;n\equiv 0(mod
m)}(-1)^{\sigma(n)},
$$

where $\sigma(n)$ is the number of 1's in the binary expansion of
$n$.

     Let $D$ be subset of primes $p$ for which $\exists x_p$ such
that $S_p(x)>0$ for $x\geq x_p$. M.Drmota and M.Skalba \cite{2}
proved that only primes $p\in(1,1000)\cap D$ are
$3,5,17,43,257,683$. In \cite{4} we conjectured that if
$p\in\mathbb{P}\backslash D$ then

\begin{equation}\label{1}
S_p(2^p)= \pm p.
\end{equation}

Below we prove this conjecture in the case of 2 is a primitive root
of $p$ (Theorem 1) and give a sufficient condition when (\ref{1})
satisfies (Theorem 2).

\section{Theorems}

\begin{theorem}\label{t1}

If $2$ is a primitive root of $p$ then $S_p(2^p)=p$.
\end{theorem}

\slshape Proof.\upshape    It is known (sf \cite{2}, formula (8) for
$y=-1$ ) that

\begin{equation}\label{2}
S_p(2^k)= \frac 1 p
\sum^{p-1}_{l=1}\prod^{k-1}_{j=0}(1-\omega_p^{l2^j}),
\end{equation}

where $\omega_p\neq 1$ is a primitive root of $1$ of the power $p$.
Since $2$ is a primitive root of $p$ then $2^j,\;j=0,1,\ldots,p-1$,
runs values $1,2,\ldots,p-1,1$ modulo $p$.
\newpage

Therefore, by (\ref{2}) we have

\begin{equation}\label{3}
S_p(2^p)= \frac 1 p
\sum^{p-1}_{l=1}(1-\omega^l_p)\prod^{p-1}_{j=0}(1-\omega_p^{lj}).
\end{equation}

Nevertheless,

$$
\prod^{p-1}_{j=1}(1-\omega_p^{lj})=\prod^{p-1}_{j=1}(1-\omega_p^j),\;\;l\in[1,p-1],
$$

and

$$
\sum^{p-1}_{l=1}\omega^l_p=-1+\sum^{p-1}_{l=0}\omega^l_p=-1.
$$
Therefore, by (\ref{3})

$$
S_p(2^p)= \frac 1 p \prod^{p-1}_{j=1}(1-\omega_p^{lj})(p-(-1))=
\prod^{p-1}_{j=1}(1-\omega_p^j).
$$

Finally, since

$$
x^p-1=(x-1)\prod^{p-1}_{j=1}(x-\omega_p^j)
$$

then

$$
\prod^{p-1}_{j=1}(x-\omega_p^j)=1+x+\ldots+x^{p-1}
$$

and thus,

$$
\prod^{p-1}_{j=1}(1-\omega_p^j)=p. \;\;\blacksquare
$$

In \cite{1} and \cite{3} it was proved by different ways that, in
particular,

$$
0< S_3(n) \leq n^{\frac{\ln{3}}{2\ln{2}}}.
$$

Using methods of \cite{3} and relations \cite{4} it could be proved
that

$$
0< S_5(n) \leq n^{\frac{\ln{5}}{4\ln{2}}}.
$$
\newpage

and

$$
|S_7(n)|\leq n^{\frac{\ln{7}}{6\ln{2}}}.
$$

In the connection with these estimates it is actual the following
theorem.

\begin{theorem}\label{t2}

If for a prime $p\geq 3$,

\begin{equation}\label{4}
|S_p(n)|\leq n^{\frac{\ln{p}}{(p-1)\ln{2}}},\;\;n\geq 2^p
\end{equation}

and $p | S(2^p)$ then $S_p(2^p)=\pm p$.
\end{theorem}

\slshape Proof. \upshape     We have

\begin{equation}\label{5}
|S_p(2^p)|< 2^{\frac{p\ln{p}}{(p-1)\ln{2}}}=p \cdot
p^{\frac{1}{p-1}}.
\end{equation}

Since $n< 2^{n-1},\;\;n\geq 3$, then by (\ref{5})

$$
|S_p(2^p)|< 2p
$$

and consequently $S_p(2^p)\in \{-p,0,p\}$.

Nevertheless, the maximal  multiple of $p$ in $[0,2^p)$ is $2^p-2$.
Thus, the number of the multiples of $p$ in $[0,2^p)$ is odd:
$\frac{2(2^{p-1}-1)}{p}+1$.

Therefore, $S_p(2^p)\neq 0.\;\;\blacksquare$

     Note that, in \cite{4} we conjectured that always $p|S_p(2^p)$.
Moreover, at least, in the case when $p$ has a primitive root $2$ we
now conjecture that together with $S_p(2^p)=p$ the relation $S_p(2^p
x)=pS_p(2x)$ and the estimate (\ref{4}) satisfy as well.

\end{document}